\documentclass{amsart}
\usepackage{righttag}
\usepackage{epsf}
\def\hepsffile{\leavevmode\epsffile}
\begin{document} 

\theoremstyle{plain}
\newtheorem{thm}{Theorem}[subsection]
\newtheorem{cor}[thm]{Corollary}
\newtheorem{lem}[thm]{Lemma}
\newtheorem{conjec}[thm]{Conjecture}
\newtheorem{prop}[thm]{Proposition}

\theoremstyle{remark}
\newtheorem{rem}{Remark}\renewcommand{\therem}{}

\theoremstyle{definition}
\newtheorem{defin}[thm]{Definition}
\newtheorem{defins}[thm]{Definitions}
\newtheorem{emf}[thm]{}

\def\lk{\protect\operatorname{lk}}
\def\C{{\mathbb C}}
\def\Z{{\mathbb Z}}
\def\Q{{\mathbb Q}}
\def\R{{\mathbb R}}
\def\N{{\mathbb N}}
\def\rpp{{\R P^2}}
\def\rpn{{\R P^n}}
\def\cpp{{\C P^2}}
\def\ca{{\C A}}
\def\ra{{\R A}}
\def\1{\hbox{\rm\rlap {1}\hskip.03in{\textrm I}}}
\def\Bbbone{{\rm1\mathchoice{\kern-0.25em}{\kern-0.25em}
        {\kern-0.2em}{\kern-0.2em}I}}
\def\p{\partial}

\font\eightrm=cmr8
\font\tenrm=cmr10
\font\eightit=cmti8
\centerline{{\em This article has been published in:\/}}
\centerline{{\em J. Knot Theory Ramifications {\bf 7} (1998), no.2, pp.
257-266.\/}}

\author[V.~Tchernov]{Vladimir Tchernov}
\address{DMATH G-66.4,
Eidgen\"osische Technische Hochschulle, CH-8092 Z\"urich, Switzerland}
\email{chernov@math.ethz.ch}
\keywords{}
\title[The Most Refined Vassiliev Invariant of Degree One]
{The Most Refined Vassiliev Invariant of Degree One\\ of Knots and
Links\\ in $\R^1$-fibrations over a surface} 

\begin{abstract}
As it is well-known, all Vassiliev invariants of degree one of a knot
$K\subset \R^3$ are trivial. There are nontrivial Vassiliev invariants of
degree one, when the ambient space is not $\R^3$.
Recently, T.~Fiedler introduced such
invariants of a knot in an $\R^1$-fibration over a surface $F$.
They take values in the free $\Z$-module generated by all the free homotopy
classes of loops in $F$. Here, we generalize them to the
most refined Vassiliev invariant of
degree one. The ranges of values of all these invariants are explicitly
described.

We also construct a similar invariant of a two-component
link in an $\R^1$-fibration. It generalizes the linking number.
\end{abstract}

\maketitle

Most proofs in this paper are postponed till the last section.

Everywhere in this text $\R^1$-fibration means a locally-trivial fibration
with fibers, homeomorphic to $\R^1$. 

We work in the differential category.

\section{Invariants of knots and links}

\subsection{Basic definitions}
We say, that a one-dimensional submanifold $L$ of a total space $N^3$ of 
a fibration 
$p:N^3\rightarrow M^2$ is {\em generic
with respect to $p$\/ }, if $p\big|_L$ is a generic immersion.
An immersion of a one-manifold into a surface is said to
be {\em generic\/}, if it has neither self-intersection points of
multiplicity greater than two, nor self-tangency  points, and at
each double point its branches are transversal to each other. An
immersion of (a circle) $S^1$ to a surface is called a {\em curve\/}.

Let $F$ be a connected smooth two-dimensional surface 
(not necessarily
compact or orientable) and $p:E\rightarrow F$ be an $\R^1$-fibration 
with oriented total space $E$.
Let $K\subset E$ be a (smooth) oriented knot, in general position with
respect to $p$.
 
\begin{defin}[Fiedler \cite{Fiedler}]\label{writhe-def}
Let $q$ be a
double point of $p(K)$. Fix an orientation on the fiber
$E_q=p^{-1}(q)$. This determines, which of the two branches of $K$, 
intersecting $E_q$, is over-crossing  and which is under-crossing. 
Define local writhe $\omega (q)$ to be one if the three-frame (under-crossing, 
over-crossing, fiber $E_q$)
agrees with the orientation on $E$ and minus one, otherwise.
(It is easy to check, that this definition does not depend on the choice of
an orientation
on $E_q$.)
\end{defin}

\subsection{Direct generalization of Fiedler's invariants}

In~\cite{Fiedler} T.~Fiedler introduced invariants of a knot $K$ in an
oriented total space of an $\R^1$-fibration $p:E\rightarrow F$. As it
follows from~\cite{Tchernov}, these invariants 
can be expressed through an invariant $U_K$,
introduced below. If $F$ is oriented, then $U_K$ also 
can
be expressed through Fiedler's invariants.  
The formulas, expressing them
through each other (see~\cite{Tchernov}), involve the values of all these
invariants on some fixed knot homotopic to $K$. 

\begin{figure}[htb]

 \begin{center}
  \epsfxsize 6cm
  \hepsffile{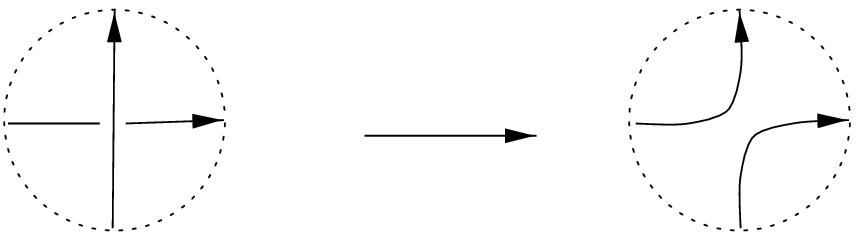}
 \end{center}
\caption{}\label{aicardi10.fig}
\end{figure}

Let $q\in p(K)$ be a crossing point. Split the curve $p(K)$ at $q$ 
according to the orientation and obtain two oriented loops
on $F$ (see Figure~\ref{aicardi10.fig}). 

\begin{defin}\label{symdef} For a 
crossing point $q$ of $p(K)$ denote by $\xi _1(q)$ and
$\xi _2(q)$ the free homotopy classes of the two loops, 
created by splitting at $q$. Let $H$ be the free $\Z$-module generated by
the set of all the free homotopy classes of oriented loops on $F$.
Define $U_K\in H$ by the following formula, 
where the summation is taken  over all the crossings, 
such that none of the two loops, created by splitting, is homotopic to a
trivial loop.  
\begin{equation}\label{eqUK}
U_K=\sum_{\{q
\in Q|\xi _1(q),\xi _2(q)\neq e\} } \omega (q)\Bigl(\xi
_1(q)+\xi _2(q)\Bigr)
\end{equation}
\end{defin} 

\begin{thm}
$U_K$ is an isotopy invariant of the knot $K$.
\end{thm}

The proof is straightforward. One checks, that $U_K$ does not
change
under all the oriented versions of the three Reidemeister moves.

\begin{emf}\label{r1homvers}
Similarly to~\cite{Fiedler}, 
one can introduce a version of $U_K$, which takes values in $\Z[H_1(F)]$. 
To obtain it, one substitutes $\xi_1(q)$ and $\xi_2(q)$
in~\eqref{eqUK} by the homology classes, realized by the corresponding loops. 
The summation should
be made over the set of all the double points of $p(K)$, such that none of the
two loops created by the splitting is homologous to $0$. 

\end{emf}

\begin{emf}
Let $p:E\rightarrow F$ be an $\R^1$-fibration over a surface. Let $K\subset
E$ be a knot generic with respect to $p$ and $q$ be a crossing point of
$p(K)$. The modification of pushing of one branch of $K$ through the other 
along a fiber $E_q$ is called the {\em modification (of the knot) along the
fiber\/} $E_q$.
\end{emf}

\begin{thm}\label{vasilyev1}(Cf. Fiedler~\cite{Fiedler}) Let $q$ be a crossing point of $p(K)$.
Denote by $i$ and $j$ the free homotopy classes 
of the two loops, created by splitting of 
$p(K)$ at $q$ according to the orientation. Under the modification along 
$E_q$ the jump of $U_K$  is
\begin{equation}\label{type1}
\begin{cases}
\pm 2\Bigl( i+j\Bigr),&\text{ if }i,j\neq e,\\
0,&\text{ otherwise. }
\end{cases}
\end{equation}
Here the sign depends on $\omega(q)$.

\end{thm}

The proof is straightforward.

\begin{cor}\label{deg1} $U_K$ is a Vassiliev invariant of degree one.
\end{cor}

To get the proof, one notices, that the first derivative of $U_K$
depends only on the free homotopy classes of the two loops, that
appear,
if one splits the singular knot (with one transverse double point) at the
double point according to the orientation. Hence, the second derivative of 
$U_K$ is identically $0$.

\subsection{The most refined Vassiliev invariant of degree one.}

\begin{emf} Unfortunately $U_K$ appears to be not the most refined Vassiliev
invariant of degree one of a knot in an $\R^1$-fibration. To show this, we 
construct two knots $K_1$ and $K_2$ and a first degree Vassiliev invariant 
$\tilde U_K$, such that $U_{K_1}=U_{K_2}$, and $\tilde U_{K_1}\neq \tilde
U_{K_2}$. 
\end{emf}

\begin{defin}[of $\tilde U_K$]\label{graphdef}
Let $\Gamma$ be an oriented figure eight graph (bouquet of two circles), $V_{\Gamma}$ be its vertex
and $E^1_{\Gamma}$ and $E^2_{\Gamma}$ be its edges.
Set $S$ to be a set of free homotopy classes of mappings of $\Gamma$ into $F$, 
factorized by an orientation preserving involution of $\Gamma$. 
Let $G$ be the free $\Z$-module generated
by $S$.
For a double point $q$ of $p(K)$ put $G_q\in S$ to be the class of the
mapping of $\Gamma$, which sends $V_{\Gamma}$ to $q$,
$E^1_{\Gamma}\cup E^2_{\Gamma}$ onto $p(K)$, according to the
orientations of the edges, and is injective on the complement of the
preimages of the double points of $p(K)$. 
Let
$S'\subset S$ be those classes, for which none of the two loops of the figure
eight graph is homotopic to a trivial loop.   
Define $\tilde U_K\in G$ by the 
following formula, where the summation is taken over the set of all the
crossings $q$ of $p(K)$, such that $G_q\in S'$.
$$\tilde U_K=\sum_{\{q\in p(K)|G_q\in S'\}}\omega (q)G_q$$
\end{defin}

\begin{figure}[htb]

 \begin{center}
  \epsfxsize 5cm
  \hepsffile{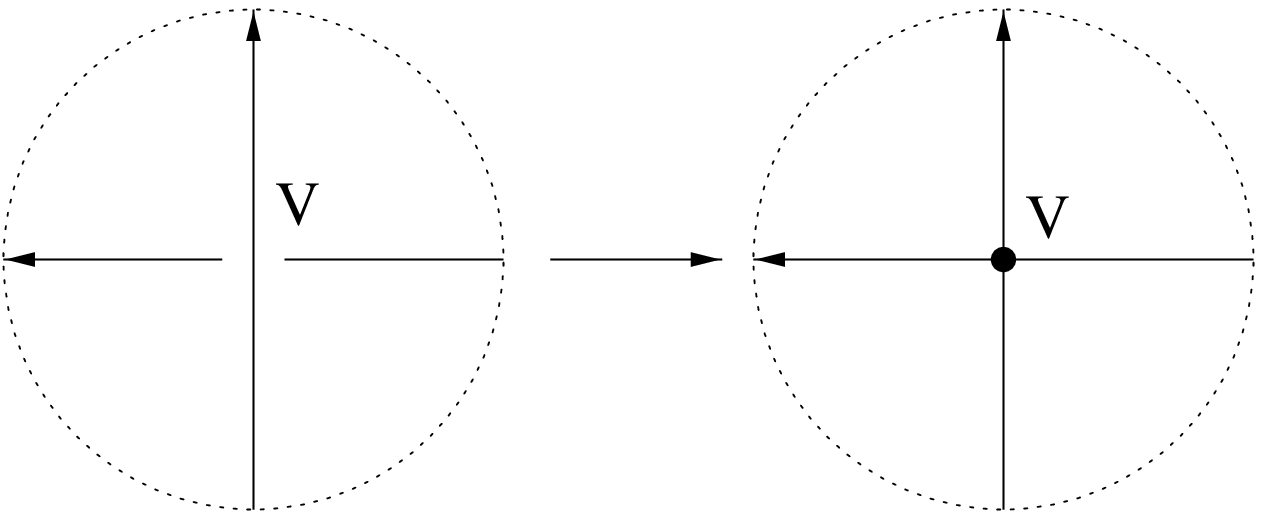}
 \end{center}
\caption{}\label{graph1.fig}
\end{figure}

Similarly to~\ref{deg1} one checks, that $\tilde U_K$ is a Vassiliev
invariant of degree one.

\begin{figure}[htb]

 \begin{center}
  \epsfxsize 7cm
  \hepsffile{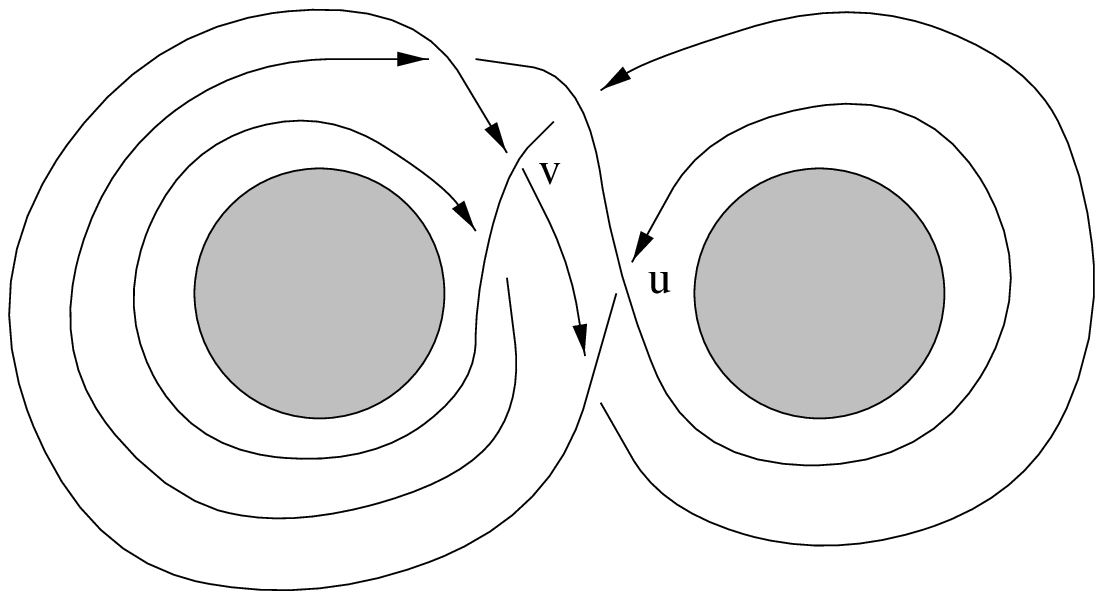}
 \end{center}
\caption{}\label{graph2.fig}
\end{figure}

Let $F$ be a disc with two holes.
Let $K_1$ be the knot, shown on Figure~\ref{graph2.fig}, and $K_2$ be the
knot obtained from $K_1$ by modifications along fibers over 
the crossing points $u$ and $v$. 
(The two shaded
discs on Figure~\ref{graph2.fig} are the two holes.) One can easily check,
that $U_{K_1}=U_{K_2}$, but $\tilde U_{K_1}\neq \tilde U_{K_2}$.

The following theorem shows, that $\tilde U_K$ invariant is the most refined 
Vassiliev invariant of degree one.

\begin{thm}\label{refine}
Let $v_1(K)$ be any Vassiliev invariant of degree one. It induces a mapping
$v_1^*:G\rightarrow \Z$, which maps a class of the projection of a singular 
knot $K'$ to $v_1(K')$. Fix some knot $K_f$. Then for any knot $K$, which is
free homotopic to $K_f$ 

\begin{equation}\label{refeq1}
v_1(K)=v_1(K_f)+\frac{1}{2}v_1^*(\tilde U_K-\tilde U_{K_f}).
\end{equation}
\end{thm}

\begin{emf}{\em Proof of Theorem~\ref{refine}.}\label{pfrefine}

One can obtain $K$ from $K_f$ by a sequence of isotopies and modifications
along fibers. Both $v_1$ and $\tilde U_K$ are invariant under isotopy. If under a
modification along a fiber $\tilde U_K$ jumps by $2G_q$, then $v_1$ jumps by 
$2v_1^*(G_q)$. (Clearly $v_1$ does not jump under modification along a
fiber,
for which one of the two loops of $G_q$ is homotopic to a trivial loop.)
The total jump of $\tilde U_K$ under the homotopy is 
$\tilde U_{K_f}-\tilde U_K$. Thus the corresponding jump of $v_1$ invariant
is $v_1(K_f)-v_1(K)=\frac{1}{2}v_1^*(\tilde U_{K_f}-\tilde U_K)$ and we
proved the theorem.
\qed
\end{emf}

It is natural to take the simplest knot in the corresponding class as
the $K_f$ knot. Unfortunately, there is no canonical way to choose one. 

As a corollary of Theorem~\ref{refine} we get, that 
for any Vassiliev invariant of degree one --- $v_1$ and 
two homotopic knots $K_1$ and $K_2$, equality 
$\tilde U_{K_1}=\tilde U_{K_2}$ 
implies $v_1(K_1)=v_1(K_2)$. 

The following theorem, characterizes the range of values $\tilde U_K$.

\begin{thm}\label{realization3}
For a singular knot $K_s$ (whose only singularity is a transverse double
point) denote by $\bar K_s$ the free homotopy class of knots, that contains 
$K_s$. For a knot $K$ denote by $G_K$ the
submodule of $G$ generated by the classes of the projections of singular 
knots $K_s$, such that $K\in \bar K_s$.

\textrm{I:} Let $K$ and $K'$ be two oriented knots, representing 
the same free homotopy class. Then $\tilde U_{K}$ and 
$\tilde U_{K'}$ are congruent modulo the $2G_{K}$ submodule.

\textrm{II:} Let $K$ be an oriented knot, $\tilde U$ be an element of $G$, 
such that it is congruent to $\tilde U_{K}$ modulo the $2G_{K}$
submodule. Then there exists an oriented knot $K'$, such that:

a) $K$ and $K'$ represent the same free homotopy class.

b) $\tilde U_{K'}=\tilde U$.
\end{thm}

For the proof of Theorem~\ref{realization3} see
Section~\ref{pfrealization3}.

\begin{emf}\label{realization1}
There is a natural mapping $\phi:G\rightarrow H$, which maps $g\in G$ 
to a formal 
sum of the free homotopy classes of the two loops of $g$.  Clearly, 
$\phi(\tilde U_K)=U_K$. (The $\ker(\phi)$ is
nontrivial and this is the reason, why $U_K$ is not the most refined
invariant of degree one.)
Using $\phi$ and Theorem~\ref{realization3} we
obtain the following characterization of the range of values of
$U_K$.

\textrm{I:} If $K$ and $K'$ are two oriented knots representing 
the same free homotopy class, then
$U_{K}$ and $U_{K'}$ are congruent modulo the $\phi(2G_K)$ submodule.

\textrm{II:} Let $K$ be an oriented knot, $U$ be an element of $H$,
 such that it is congruent to  $U_{K}$ modulo the $\phi(2G_K)$ submodule. 
Then, there exists an oriented knot
$K'$,
such that:

a) $K$ and $K'$ represent the same free homotopy class.

b) $U_{K'}=U$.
\end{emf}

\subsection{Partial linking polynomial}\label{aic-def}
Let $\Theta$ be an annulus.
Consider a solid torus $T$ embedded into $\R^3$, and a projection 
$p :\R^3\rightarrow \R^2$, such that $\Im(p \big|_T)$ is
homeomorphic to $\Theta$. Let $K \subset T$ be an oriented
knot, in general position with respect to $p$.
We denote by $i_1(q)\text{ and }i_2(q)$ the homology classes
in $H_1(\Theta)$ of the two loops, that are created by splitting of $p (K)$
at the double point $q$. Since $H_1(\Theta)=\Z$ we can
consider $i_1(q)$ and $i_2(q)$ as integer numbers.
\begin{defin}[Aicardi \cite{Aicardi}] Set partial
linking polynomial $A(K)$ (originally in~\cite{Aicardi} it was denoted by
$s[K]$) to be a finite Laurent polynomial, defined 
by the following formula $$
 A(K)=\sum_{\{q\in
Q|i_1(q),i_2(q)\ne 0\} }\frac{1}{2}\Bigl(\omega
(q)(t^{i_1(q)}+t^{i_2(q)})\Bigr).$$
Below by $a_i$ we denote the coefficient of $t^i$ in $A(K)$.
\end{defin}

\begin{emf}\label{relation}
The set of all the free homotopy classes of oriented loops in $\Theta$
coincides with $H_1(\Theta)$.
One can easily see, that  $U_K$ is mapped to $2A(K)$ under the natural
isomorphism $\psi: H\rightarrow \Z[q,q^{-1}]$. 

The fact, that 
$\pi_1(T)=\Z$ allows one to reconstruct an element $g\in G$
from the homology classes of the two loops of it. Thus, 
in this case $\tilde U_K$
invariant can also be reconstructed from $A(K)$. 
\end{emf}

\begin{emf}[Aicardi \cite{Aicardi}]  Let $h\in \Z$ be the image of
$[p(K)]$ (the homology class realized by $p(K)$)  
under the natural identification of $H_1(\Theta)$ with $\Z$. 
Then $a_0=a_h=0$ and $a_i=a_{h-i}$ for an arbitrary $i\in \Z$.
\end{emf}

\begin{emf}
One can see, that the very definition of $A(K)$ depends on the embedding
of $T$ into $\R^3$. It is well known, that the group of 
orientation preserving
autohomeomorphisms of $T$, factorized by isotopy relation,
is isomorphic to $\Z$. It is generated by the class of an autohomeomorphism
$\Phi$, 
that extends a positive Dehn twist along a meridian of
$\p T$. That is cutting $T$ along a meridional disc, twisting
by $2\pi $ in a positive direction and gluing back. Replacement of the 
embedding of $T$ to $\R^3$ by an isotopic one does not change $A(K)$. 
Embeddings of all isotopic classes can be obtained from the given one by a 
composition with $\Phi ^n$ for some $n\in \Z$.
\end{emf}
 
Let $A'(K)$ be the partial linking polynomial calculated, 
after we compose our embedding of $T$ with $\Phi$. 
Put

$$\Delta A(K)=A'(K)-A(K).$$

Let $h\in \Z$ be the homology class realized by $p(K)$.

\begin{thm}\label{twist}

$$\Delta A(K)=\begin{cases} -|h|(t^1+t^2+\dots +t^{h-1}), &
 \text{if }h>0\\ 
-|h|(t^{-1}+t^{-2}+\dots +t^{h+1}), &\text{if } h<0\\
0, &\text{if } h=0\end{cases}$$

\end{thm}

For the proof of Theorem~\ref{twist} see Section~\ref{pftwist}.

As we can make the composition of our embedding with $\Phi ^n$, for any
$n\in \Z$, we obtain the following.

\begin{emf}$A(K)$ as an invariant of the topological pair $K \subset
T$ is defined up to an addition of $\Delta A(K)$. 
Thus, an $A(K)$ invariant of a knot $K$, could be said to be in a canonical
form, if it satisfies the following conditions:
$$\begin{cases} 0\leq a_1<h & \text{for }h>0, \\ 
0\leq a_{-1}< |h| &\text{for }h<0, \\  
\end{cases}$$
If $h=0$, then $A(K)$ is always in the canonical form.

\end{emf}

\begin{thm}\label{character}
Fix $h\in \Z$. Let $P_h$ be a subset of all finite Laurent polynomials
$\sum _{i=i_1}^{i_2} p_it^i$, satisfying the following properties:

a) $p_0=p_h=0$ 

b) $\forall j\in \Z\quad  p_j=p_{h-j}$ 

c) if $h=2k$ for some $k \in \Z$ then $p_k$ is odd.

Then $P_h$ is the range of values of the partial linking
polynomial for knots homologous to $h$.
\end{thm}

For the proof of Theorem~\ref{character} see Section~\ref{pfcharacter}.

\subsection{Invariant of links}

\begin{defin}[of $U_L$]
Let $p:E\rightarrow F$ be an $\R^1$-fibration, of an oriented space $E$ over
a surface.
Let $\Gamma$ be an oriented figure eight graph (bouquet of two circles), $V_{\Gamma}$ be its vertex
and $E^1_{\Gamma}$ and $E^2_{\Gamma}$ be its edges.
Set $\bar S$ to be a set of all the free homotopy classes of 
mappings of $\Gamma$ into $F$. 
Denote by $\bar G$ the free $\Z$-module generated
by $\bar S$. 
Let $K_1\cup K_2=L\subset E$ be an oriented two-component link, in
general position with respect to $p$.
Note, that local writhe $\omega (q)$ 
is well
defined for a point $q\in
p(K_1)\cap p(K_2)$.
Let $\bar G_q\in \bar S$ be the class of the mapping of $\Gamma$ 
onto $p(K_1)\cup p(K_2)$, which maps $V_{\Gamma}$ to $q$, 
$E^1_{\Gamma}$ to $p(K_1)$, $E^2_{\Gamma}$ to $p(K_2)$ 
(according to the orientations of the edges) and is injective on the
complement of the preimage of the double points of $p(L)$. 
Define 
$U_L\in \bar G$ by the
following formula, where the summation is taken over $p(K_1)\cap p (K_2)$ 
\begin{equation}\label{eqUKij}
U_L=\sum_{q\in p(K_i)\cap p(K_j)}\omega(q) \bar G_q
\end{equation}
\end{defin}

\begin{thm}\label{correctUL}
$U_L$ is an isotopy invariant of the link $L$.
\end{thm}

The proof of Theorem~\ref{correctUL} is straightforward. One just has to
check, that $U_L$ is invariant under all the oriented versions
of the Reidemeister moves.

\begin{emf}
If $E=\R^3$ and $F=\R^2$, then $\bar G=\Z$ (as $\pi_1(\R^2)=e)$. Under this
identification $U_L=2\lk(K_1,K_2)$, where $\lk(K_1,K_2)$ is the linking
number of the two knots.
\end{emf}

\begin{emf}
Let $L=K_1\cup\dots\cup K_n\subset E$ be a generic 
$n$-component oriented link. For $i>j$ ($i,j\in\{1,\dots,n\}$) set $L_{ij}$
to be the two component sublink of $L$, consisting of $K_i$ and $K_j$.
Similarly to Theorem~\ref{refine}, one can see, that the 
ordered set of the invariants $U_{K_i}$ and $U_{L_{ij}}$ ($i>j$) is the most
refined degree one Vassiliev invariant of $L$.
\end{emf}

\section{Proofs}

\subsection{Proof of Theorem~\ref{realization3}.}\label{pfrealization3}
\textrm{I:} $K'$ can be obtained from $K$ by a sequence of isotopies and
modifications along fibers. Isotopies do not change $\tilde U$. 
The modifications 
change $\tilde U$ by
elements of $2G_{K}$. Thus, the first part of the theorem is proved.

\begin{figure}[htbp]
 \begin{center}
  \epsfxsize 3cm
  \hepsffile{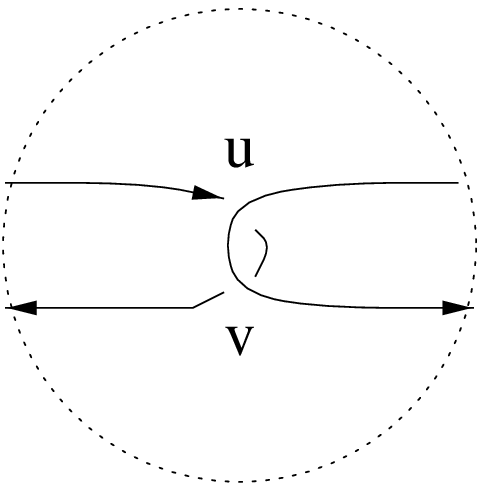}
 \end{center}
\caption{}\label{bite.fig} 
\end{figure}

\textrm{II:} We prove that for any $g\in G_{K}$ there exist
two knots $K_1$ and $K_2$ such, that they represent the same free homotopy
class as $K$ and 
$$\tilde U_{K_1}=\tilde U_K-2g$$ $$\tilde U_{K_2}=\tilde U_K+2g$$
Clearly, this implies the second statement of the theorem.
To obtain the two knots we isotopically deform $K$ so that $\pi(K)$ bites
itself in the projection (as it is shown in Figure~\ref{bite.fig}) and
$G_u=G_v=g$. To obtain $K_1$, one performs a fiber modification along 
$\pi^{-1}(u)$. To obtain $K_2$, one performs a fiber modification along 
$\pi^{-1}(v)$.
This finishes the proof of Theorem~\ref{realization3}.
\qed

\subsection{Proof of Theorem~\ref{twist}.}\label{pftwist}
\begin{figure}[htb]
 \begin{center}
  \epsfxsize 6cm
  \hepsffile{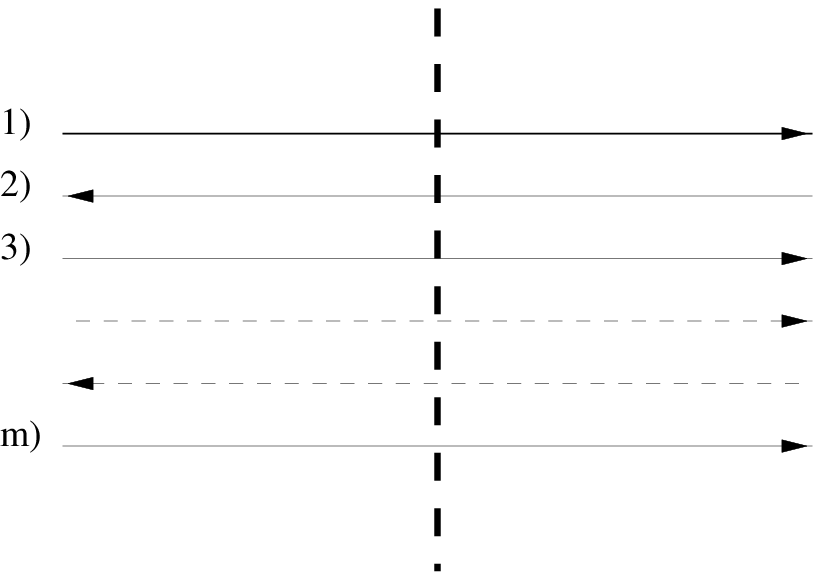}
 \end{center}
\caption{}\label{aicardi4.fig}
\end{figure}

\begin{figure}[htb]
 \begin{center}
  \epsfxsize 11cm
  \hepsffile{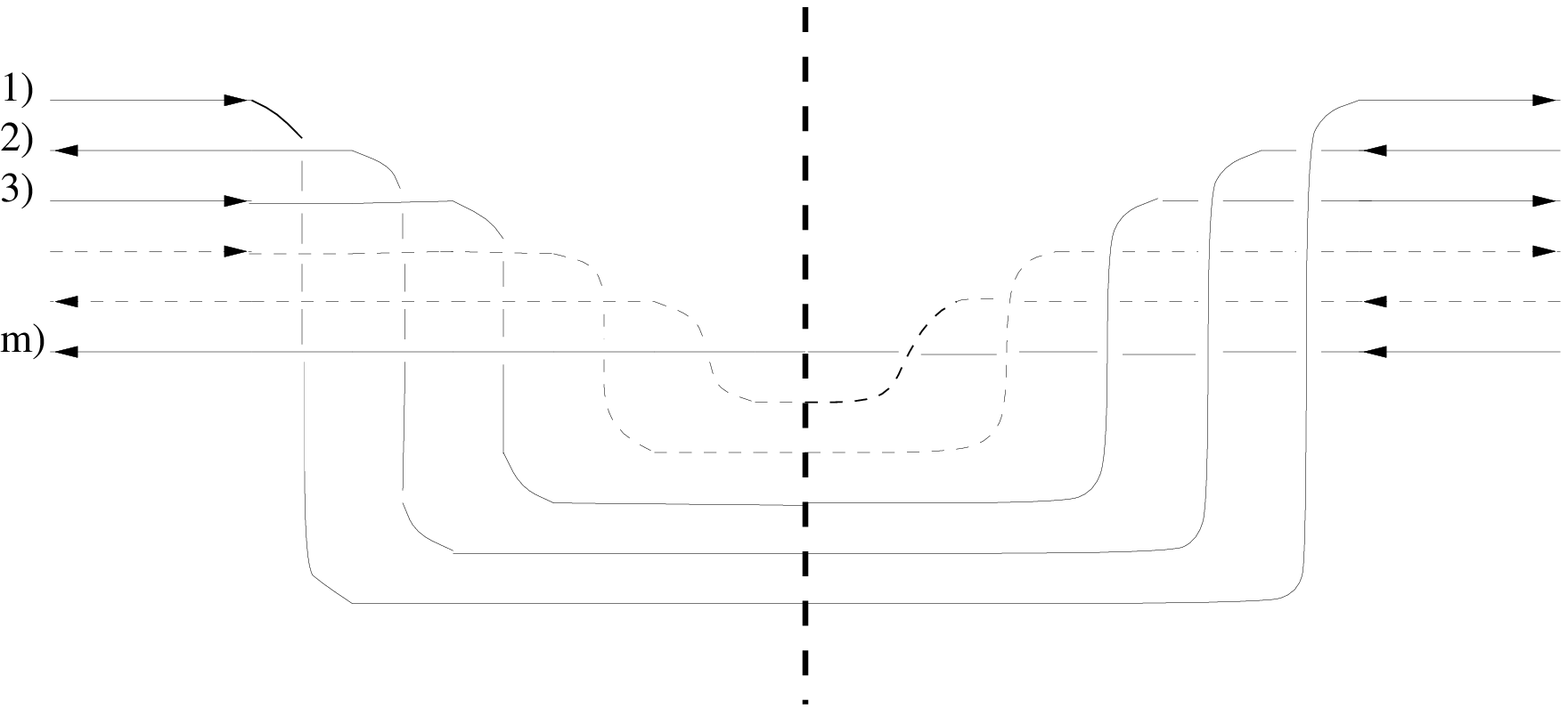}
 \end{center}
\caption{}\label{aicardi2.fig}
\end{figure}

Let $D$ be a meridional disc along the boundary of which, 
we performed the positive Dehn twist (used to define $\Phi$). 
Assume, that all the branches of K, which cross $D$, are perpendicular
to it and are located on different levels (see Figure~\ref{aicardi4.fig}). 
Using second Reidemeister moves 
transform the diagram in such a way, that if we 
traverse $K$ along the orientation, 
then the branches cross $D$ in the order shown in Figure~\ref{aicardi4.fig}. 
(The thick dashed line in Figure~\ref{aicardi4.fig} is $p(D)$).  
After we compose the embedding of $T$ with $\Phi$, 
the diagram  will be changed, as it is shown in Figure~\ref{aicardi2.fig}.

Note, that under the modification of pushing of one branch of the knot
through the other, which  happens outside of the neighborhood of $D$ (shown in
Figure~\ref{aicardi2.fig}) $A(K)$ and $A'(K)$ change in the same way. Hence,
their difference is preserved. Thus, we can assume that our knot $K$ 
has an ascending diagram. 
After a simple calculation we get the desired result.
\qed

\subsection{Proof of Theorem~\ref{character}.}\label{pfcharacter}
 The relation between $U$ and $A$ invariants, shown in~\ref{relation}, 
allows one to use~\ref{realization1} in the case of 
a partial linking polynomial. There is a natural 
bijection between one-dimensional homology classes of $T$ 
and free homotopy classes of oriented loops in $T$. 

Thus we get, that:

a) If $K$ and $K'$, are such that $[p(K)]=[p(K')]=h$, then 
$A(K')$ and $A(K)$ are congruent modulo the
additive subgroup generated by all the elements of type
 \begin{equation}\label{polchange}
\pm(t^j+t^{h-j}) \text{     for } j\not\in\{ h,0\}
 \end{equation}
(Note that if $h=2j$ then this expression is equal to $\pm 2t^j$.)    

b) Let $K$ be a knot (with $[p(K)]=h$), and let $A$ be a finite Laurent
polynomial congruent to $A(K)$ modulo the additive subgroup,
generated by all the elements of type~\eqref{polchange}. Then there exists a 
knot $K'$, such that $[p(K')]=h$ and $A(K')=A$.

\begin{figure}[h]
 \begin{center}
  \epsfxsize 4cm
  \hepsffile{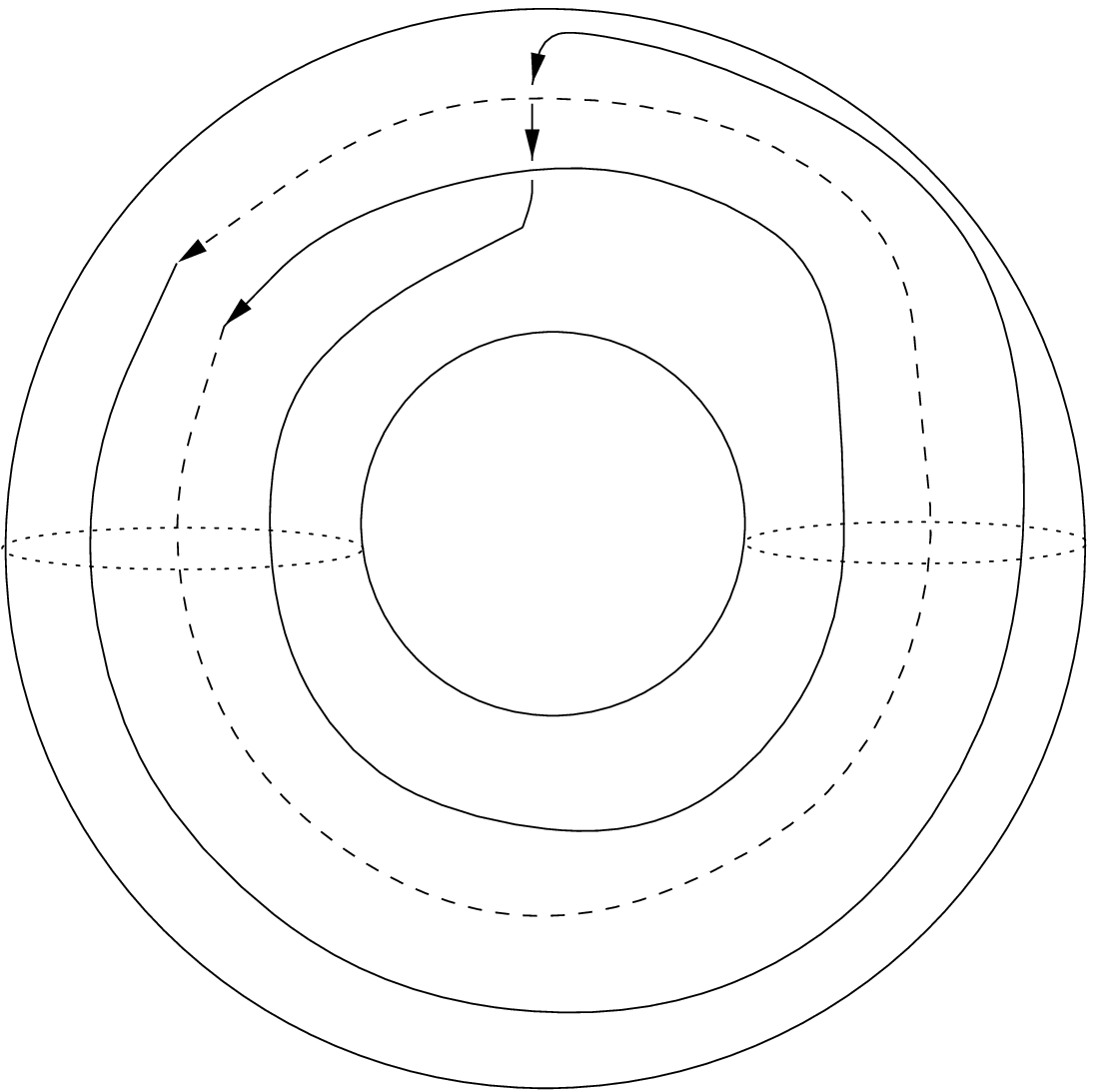}
 \end{center}
\caption{}\label{aicardi6.fig}
\end{figure}

Thus, if $K$ and $K'$ are knots such, that $[p(K)]=[p(K')]=h$ and 
$A(K)\in P_h$, then $A(K')\in P_h$. 
And vice versa, if for some $p_h\in P_h$ there exists a knot $K_{p_h}$, such
that $[p(K_{p_h})]=h$ and $A(K_{p_h})=p_h$, then such a knot exists for any $\tilde
 p_h\in P_h$. Hence, to prove the theorem it is sufficient to show, that 
for any $h\in \Z$ there exists a knot $K_h$, such that $[p(K_h)]=h$ and
$A(K_h)\in P_h$. Let $K_h$ be a knot, that rotates $h$ times in $T$ and has 
an ascending diagram (see Figure~\ref{aicardi6.fig}). The $A$ invariant of
it is equal to~\eqref{simplepol} and it belongs to $P_h$.
 \begin{equation}\label{simplepol} 
   \begin{cases} t^1+t^2+\dots +t^{h-1}, &\text{if } h>0 \\
                 t^{-1}+t^{-2}+\dots +t^{h+1}, &\text{if } h<0\\
		 0, &\text{if } h=0\\
    \end{cases}
 \end{equation}   
This finishes the proof of Theorem~\ref{character}.
\qed  

\noindent{\bf Acknowledgements}
\centerline{}

I am deeply grateful to Oleg Viro for the inspiration of this work and all the
enlightening discussions. I am thankful to Francesca Aicardi, Thomas
Fiedler and Michael Polyak for all the valuable discussions we had.


\begin{thebibliography}{99999}

\bibitem{Fiedler}
T.~Fiedler, {\em A small state sum for knots}, Topology {\bf 30} (1993)
no.2, 281-294
\bibitem{Tchernov}
V.~Tchernov {\em First degree Vassiliev invariants of knots in $\R^1$- and
$S^1$-fibrations} preprint, Uppsala, Sweden 1996
\bibitem{Aicardi}
F.~Aicardi, {\em Invariant Polynomial of Framed Knots in the Solid Torus and
its Applications to Wave Fronts and Legendrian Knots}, J. of Knot Th. and
Ramif, Vol 15, No. 6, (1996), 743-778  
\end{thebibliography}
\end{document}